\documentclass[11pt]{amsart}
\usepackage{geometry} 
\usepackage{graphicx,xcolor}
\usepackage[all]{xy}
\usepackage{amssymb}
\usepackage[shortlabels]{enumitem}
\usepackage{hyperref}
\usepackage{cleveref}
\usepackage{MnSymbol}
\usepackage{tikz-cd}
\usepackage{overpic}
\usepackage{caption}

\newtheorem{theorem}{Theorem}[section]
\newtheorem{proposition}[theorem]{Proposition}
\newtheorem{cor}[theorem]{Corollary}
\newtheorem{lemma}[theorem]{Lemma}
\newtheorem{remark}[theorem]{Remark}
\newtheorem{conjecture}[theorem]{Conjecture}
\newtheorem{definition}[theorem]{Definition}
\newtheorem{exmp}[theorem]{Example}

\newcommand{\crefparts}[2]{%
  \hyperref[#2]{\cref{#1}.(\ref*{#2})}%
}
\newcommand{\Crefparts}[2]{%
  \hyperref[#2]{\Cref*{#1}.(\ref*{#2})}%
}

\newcommand{\N}{\mathbb{N}}
\newcommand{\Z}{\mathbb{Z}}

\newcommand{\C}{\mathbb{C}}

\newcommand{\Pell}{\mathcal{P}}
\newcommand{\iso}{\cong}
\newcommand{\KHI}{\mathit{KHI}}
\newcommand\SU{\mathrm{SU}}

\newcommand{\bdefn}{\begin{definition}}
\newcommand{\edefn}{\end{definition}}

\usepackage{todonotes}

\author{Sudipta Ghosh}
\address{University of Notre Dame\\ USA}
\email{sghosh7@nd.edu}

\author{Zhenkun Li}
\address{Academy of Mathematics and Systems Science\\Chinese Academy of Sciences}
\email{zhenkun@amss.ac.cn}

\author{Juanita Pinz\'on-Caicedo}
\thanks{The third author was partially supported by Simons Collaboration grant 712377.}
\address{University of Notre Dame\\ USA}
\email{jpinzonc@nd.edu}

\title{$SU(2)$-representations of Branched Covers}

\begin{document}
\begin{abstract}
We study the existence of irreducible $SU(2)$-representations for cyclic branched covers of knots in $S^3$. Our main result establishes that if $K$ is a non-trivial prime knot and $d$ is an integer such that $d \geq 2$ and $\Sigma_d(K)$ is an integer homology sphere, then $\pi_1(\Sigma_d(K))$ admits an irreducible $SU(2)$-representation, whenever $K$ satisfies one of two conditions: either $K$ is $2$-periodic, or $K$ can be represented as the closure of a tangle adapted to a $d\times d$ SICUP matrix. The first condition leverages a commuting trick for covering spaces to realize higher-degree branched covers as 2-fold covers, allowing us to apply recent results of Kronheimer-Mrowka and others. The second condition uses equivariant surgery descriptions and the $\nu^\sharp$ invariant from instanton Floer homology. As applications, we provide new infinite families of hyperbolic integer homology spheres admitting irreducible representations, including examples where previously known criteria fail.
\end{abstract}
\maketitle

\section{Introduction}

The study of $SU(2)$-representations of 3-manifold groups stands at the intersection of topology, geometry, and gauge theory, providing fundamental insights into the structure of 3-manifolds. Understanding when a 3-manifold admits irreducible representations---homomorphisms from its fundamental group to $SU(2)$ with non-abelian image---has emerged as a central problem in low-dimensional topology, with deep connections to instanton Floer homology, hyperbolic geometry, and quantum invariants.

A compelling conjecture, appearing as Problem 3.105(A) on Kirby's problem list, asks whether the $SU(2)$ representation variety is nontrivial for every compact 3-manifold with nontrivial fundamental group. For integer homology 3-spheres, this specializes to:

\begin{conjecture} Let $Y$ be a closed oriented 3-manifold satisfying $H_i(Y;\Z)=H_i(S^3;\Z)$. If $Y$ is not $S^3$, then there exists a homomorphism $\rho: \pi_1(Y)\to SU(2)$ with non-abelian image.
\end{conjecture}

The significance of this conjecture extends beyond pure topology. Since generators of the instanton chain complex can be (roughly) identified with representations of 3-manifold groups, instanton theory is intimately connected to the fundamental group---an invariant that determines 3-dimensional diffeomorphism types up to finite ambiguity. This connection proved crucial in Kronheimer-Mrowka's resolution of Property P \cite{kronheimer-mrowka:P}, where they established that homology spheres obtained as Dehn surgery along non-trivial knots in $S^3$ admit irreducible $SU(2)$-representations.

Any irreducible integer homology 3-sphere decomposes along a finite collection of embedded tori into pieces admitting either Seifert fibrations or complete hyperbolic structures of finite volume \cite{en}. This decomposition suggests three natural cases for proving Conjecture 1.1:

\begin{enumerate}
\item Seifert fibered homology spheres $Y=\Sigma(a_1,a_2,\ldots,a_n)$
\item Toroidal homology spheres 
\item Hyperbolic homology spheres
\end{enumerate}

The first case was resolved in the 1980s by Fintushel-Stern \cite{fs} and Kirk-Klassen \cite{kirk-klassen}, and the second case independently by Lidman-Pinz\'on-Caicedo-Zentner \cite{toroidal} and Baldwin-Sivek \cite{bs-splicing}. The present article is a modest attempt to tackle the last case.

Our approach involves studying cyclic branched covers of knots in $S^3$. By Thurston's Orbifold Theorem, if a knot exterior $S^3\setminus K$ admits a complete hyperbolic structure, then for $n\geq 3$ the $n$-fold cyclic branched cover $\Sigma_n(K)$ typically admits a hyperbolic structure (with the exception of $\Sigma_3(4_1)$). Rather than focus only on the case of hyperbolic knots, we identify specific properties of a knot $K$ that guarantee its branched covers admit irreducible representations. However, it is worth mentioning that since $n$-fold covers branched along torus knots are Seifert fibered, and those branched along satellite knots are toroidal, the existence of representations follows from the previously mentioned results.

\begin{theorem}\label{mainthm} Let $K\subset S^3$ be a non-trivial prime knot and $d>2$ be such that $\Sigma_d(K)$ is an integer homology sphere. If either
\begin{enumerate}
\item\label{2per} $K$ has cyclic period 2, or
\item\label{nu} $K$ can be represented as the closure of a tangle adapted to a $d\times d$ SICUP matrix,
\end{enumerate}
then $\Sigma_d(K)$ admits an irreducible $SU(2)$-representation.
\end{theorem}

To demonstrate that the second condition is non-trivial, we develop an algorithm to construct an infinite family of knots adapted to any given $d\times d$ SICUP (Symmetric, Integral, Circulant, Unimodular, Positive definite) matrix. Furthermore, we fully classify all $5\times 5$ SICUP matrices, revealing an unexpected connection between topology and Pell's equation.

Our result can be counted as one in a list of results relating topological properties of a knot $K\subset S^3$ with the existence of $SU(2)$-representations of its $d$-fold cover $\Sigma_d(K)$. Other results in this direction include the non-triviality of the signature function at an $n$-th root of unity \cite{collin}, and quasipositivity and non-sliceness of $K$ \cite{bs-stein}.\\

\textbf{Outline:} The paper is organized as follows. \Cref{periodic} establishes our periodicity results, showing how 2-periodic knots lead to representations for their branched covers. \Cref{instantons} introduces the $\nu^\sharp$ invariant and proves our main technical theorem connecting SICUP matrices to representations. \Cref{sicup} provides explicit constructions of knots adapted to SICUP matrices, yielding infinite families of hyperbolic homology spheres with irreducible representations. \\

\textbf{Acknowledgments:}  We thank Tye Lidman, Mike Miller Eismeier, Steven Sivek, and Fan Ye for comments made on an earlier version of the article. 



\section{Periodic Knots and Branched Covers}\label{periodic}
The 2-periodicity condition included in the first part of our theorem is related to the existence of irreducible representations for 2-fold covers. In this section we recall these latter results and explain clearly the way periodicity plays a role. We start with a basic lemma relating the fundamental groups of a branched cover and its base. We include this proof for the sake of completeness. 

\begin{lemma}\label{surjections} Let $Y$ be an integer homology sphere, $J\subseteq Y$ be a knot, and $d\geq 1$ an integer. Denote by $\Sigma=\Sigma_d(Y,J)$ the $d$-fold cover of $Y$ branched over $J$. There exists a surjection $\pi_1\left(\Sigma_d\left(Y,J\right)\right) \to \pi_1\left (Y\right)$\end{lemma}

\begin{proof} Let $p: \Sigma \setminus N(\widetilde{J})\to Y\setminus N(J)$ denote the covering map that gives rise to the branched cover. Let $\pi=\pi_1\left(Y\setminus N(J)\right)$ and $\widetilde{\pi}=p_*\left(\pi_1\left(\Sigma\setminus N(\widetilde{J})\right)\right)$, where $\widetilde{J}$ denotes the preimage of $J$ in $\Sigma$, and $p_*$ is the induced map on fundamental groups. Let $\mu$ be a meridian of $J$ in $Y$, and denote by $\llangle\mu\rrangle$ and $\llangle\mu^d\rrangle$ the subgroups normally generated by $\mu$ and $\mu^d$, respectively. With that notation we have $$\pi_1(Y)=\pi/\llangle\mu\rrangle, \text{ and } p_*\left(\pi_1 (\Sigma)\right)=\widetilde{\pi}/\llangle\mu^d\rrangle.$$ 
To define a homomorphism $g:\pi_1 (\Sigma)\to \pi_1(Y)$, consider the chains of inclusions of normal subgroups $\llangle\mu^d\rrangle\trianglelefteq \widetilde{\pi} \trianglelefteq\pi$ and $\llangle\mu^d\rrangle\trianglelefteq \llangle\mu\rrangle\trianglelefteq\pi$. By the isomorphism theorems for groups, these give rise to the following short exact sequences of groups: 

\begin{equation}
\xymatrix@R=2ex@C=2ex{
1\ar[dr]&&&&1\\
&\widetilde{\pi}/\llangle\mu^d\rrangle\ar@{-->}[rr]^{g}\ar[dr]^{\iota} & & \pi/\llangle\mu\rrangle\ar[ur]&\\
&&\pi/\llangle\mu^d\rrangle\ar[ur]^{\varphi}\ar[dr]^{h}&&\\
&\llangle\mu\rrangle/\llangle\mu^d\rrangle\ar[ur]&&\pi/\widetilde{\pi}=\Z/d\ar[dr]&\\
1\ar[ur]&&&&1\\
}
\end{equation}
The horizontal dashed map at the top is the desired homomorphism $g:\pi_1 (\Sigma)\to \pi_1(Y)$, and is realized as the composition $\varphi\circ\iota$, where $\iota$ and $\varphi$ are induced by inclusions. The map labeled by $h$ is induced by the Hurewicz map mod $d$. To show that $g$ is surjective, let $\gamma\in\pi$ and $l=lk_Y(\gamma,J)$. Then $\varphi(\gamma\mu^{-l}\llangle\mu^d\rrangle)=\gamma\llangle\mu\rrangle$ with $\gamma\mu^{-l}\in\widetilde{\pi}$, since $h(\gamma\mu^{-l})=0$. This gives $$g\left(\gamma\mu^{-l}\llangle\mu^d\rrangle\right)=\gamma\llangle\mu\rrangle,$$ thus establishing the desired surjectivity. 

\end{proof}

\begin{theorem}
\label{2fold-thm} Let $Y$ be an integer homology sphere, and let $J\subseteq Y$ be a non-trivial knot. Denote by $\Sigma=\Sigma_2(Y,J)$ the 2-fold cover of $Y$ branched over $J$. If $\Sigma_2(Y,J)$ is an integer homology sphere, then it admits irreducible $SU(2)$-representations.
%
\end{theorem}

\begin{proof} If $Y$ itself admits an irreducible $SU(2)$-representation $\rho$, then the composition of $\rho$ with the surjection from \Cref{surjections} produces an irreducible representation for $\Sigma$. If $Y=S^3$, then the result follows from \cite[Corollary 7.17]{KM-sutures}, and \cite[Theorem 3.1]{CNS}. To establish the general case, let $Y\neq S^3$ and suppose $Y$ does not admit any irreducible $SU(2)$ representations. Assume that $\Sigma_2(Y,J)$ does not admit an irreducible $SU(2)$-representation and argue by contradiction: by \cite[Corollary 4.10]{baldwin2018stein}, we must have 
\[
\text{dim}_{\mathbb{C}}(I^{\#}(Y)) = 1.
\] 
Then \cite[Theorem 4.1]{smalldehn} implies that
\[
\dim_{\mathbb{C}} KHI(Y,J) = 1 = |H_1(Y;\mathbb{Z})|.
\]
Thus, \cite[Theorem 1.5]{li2024enhanced} implies that $J$ is the unknot, contradicting the hypothesis.

\end{proof}

\Crefparts{mainthm}{2per} is related to the problem of realizing 3-manifolds as cyclic branched coverings of different knots and different degrees. Indeed, work of Reni-Zimmermann shows that in many cases, if $M$ is a cyclic branched cover of more than one knot, then the quotient knots are symmetric \cite{reni-zim}. The specific type of symmetry we consider is the following:

\begin{definition}\label{def:2per} A knot $K\subseteq S^3$ is said to be 2-periodic if there exists a self-diffeomorphism $f:S^3\to S^3$ such that
\begin{enumerate}
\item $f^2=id$,
\item the set $A$ of fixed points of $f$ is a curve,
\item $f(K)=K$
\item $A\cap K=\emptyset$.
\end{enumerate}
In this case, $A$ is said to be a 2-periodicity axis for $K$, and $w=lk(K,A)$ is the winding number (with respect to $f$).
\end{definition}

We remark that by the resolution of the Smith conjecture \cite{smith-conj}, the fixed point set $A$ is an unknotted curve, and the quotient space of $S^3$ under the action of the group generated by $f$ is $S^3$ itself. As a consequence, the original copy of $S^3$ can be regarded as the 2-fold cover of $S^3$ (which is $S^3/\langle f\rangle$) branched over the quotient knot $A_q$ of $A$.

\begin{theorem}\label{2per-to-2fold} Let $K\subseteq S^3$ be a non-trivial knot. If $K$ is 2-periodic with periodicity axis $A$, and in addition $w=lk(K,A)$ is relatively prime to $d$, then $\Sigma_d(K)$ can be realized as a 2-fold cover of a homology sphere branched over a knot.
\end{theorem}

\begin{proof}  Let $f:S^3\to S^3$ be the involution that gives the 2-periodicity of $K$. Denote by $A$ the set of fixed points of $f$, by $q$ the quotient map $S^3\to S^3/\langle f\rangle$, and by $K_q\sqcup A_q$ the link obtained as $q(K\sqcup A)$. In order to understand branched covers, we will examine covering spaces of $X=S^3 \setminus N(K_q\sqcup A_q)$. More precisely, define the maps
\begin{align*}
\beta:&\pi\longrightarrow \Z/d & \alpha&:\pi\longrightarrow \Z/2\\
&\gamma\longrightarrow lk(\gamma,K_q) \pmod{d}& &\gamma\longrightarrow lk(\gamma,A_q) \pmod{2}.
\end{align*}
Denote by $\widetilde{K}\sqcup \widetilde{A}$ the preimage of $K\sqcup A$ in $\Sigma=\Sigma_d(K)$. Then, by the classification of covering spaces, $S^3\setminus N\left(K\sqcup A\right)$ is the covering space of $X$ corresponding to $\ker(\alpha)$, and $\Sigma\setminus N\left(\widetilde{K}\sqcup \widetilde{A}\right)$ to $\ker(\alpha)\cap\ker(\beta)$. Furthermore, if $V$ is the covering space of $X$ corresponding to $\ker(\beta)$, then $\Sigma\setminus N\left(\widetilde{K}\sqcup \widetilde{A}\right)$ is a 2-fold cyclic cover of $V$. An appropriate filling of $V$ will then accomplish the goal of realizing $\Sigma$ as a 2-fold branched cover; to understand the filling, one first has to examine the structure of $\partial V$ as a covering space of $\partial N(K_q)\sqcup\partial N(A_q)$. The assumption that $w=lk(K,A)$ is relatively prime to $d$ guarantees that the longitudes of both $K_q$ and $A_q$ map to an invertible element of $\Z/d$ under $\beta$. Thus, the boundary of $V$ is homeomorphic to two disjoint copies of $S^1\times S^1$, denoted by $\partial V_{K_q}$ and $\partial V_{A_q}$. The manifold obtained by filling $\partial V_{A_q}$ by gluing a solid torus $D^2\times S^1$ so that a meridional curve $\partial D^2\times \{1\}$ is identified with the preimage of a meridian of $A_q$ is then a covering space of $S^3\setminus N(K_q)$ corresponding to $\ker(\beta)$. If we then fill it with another solid torus in such a way that a meridional disk bounds the preimage of the product of $d$ meridians of $K_q$, we get exactly $\Sigma_d(K_q)$. Then, notice that $\Sigma$ is the 2-fold cover of $\Sigma_d(K_q)$ branched over the core of the filling of $\partial V_{A_q}$. The diagram presented in \cref{chaincovers} summarizes the chain of covering spaces used in the proof.

\begin{equation}\label{chaincovers}
\xymatrix{
&\Sigma\setminus N\left(\widetilde{K}\sqcup \widetilde{A}\right)\ar[dl]^{d-1}\ar[dr]_{2-1}&\\
S^3\setminus N\left(K\sqcup A\right)\ar[dr]^{2-1}_q&&V=\Sigma_{d}\left(K_q\right)\setminus N\left(K'\sqcup A'\right)\ar[dl]_{d-1}\\
&S^3\setminus N(K_q\sqcup A_q)&\\
}
\end{equation}
\end{proof}

\Crefparts{mainthm}{2per} follows immediately from the following statement.

\begin{cor} Let $K\subseteq S^3$ be a 2-periodic non-trivial prime knot with periodicity axis $A$. For $d\geq 3$, we consider $\Sigma=\Sigma_d(K)$, the $d$-fold cover of $S^3$ branched over $K$. If $\Sigma_d(K)$ is an integer homology sphere, then it admits an irreducible $SU(2)$ representation.
\end{cor}

\begin{proof} If $w=lk(K,A)$ is relatively prime to $d$, then by \Cref{2per-to-2fold}, $\Sigma_d(K)$ is the 2-fold cover of $\Sigma_d(K_q)$ branched over the knot $A'$ which is the preimage of $A_q$. The corollary is thus a consequence of \Cref{2fold-thm} as follows:
\begin{enumerate}
\item The case when $K_q$ is knotted and $\Sigma_d(K_q)$ admits $SU(2)$ representations is a direct consequence of \Cref{surjections}.
\item The case when $K_q$ is unknotted follows from \cite[Corollary 7.17]{KM-sutures} and \cite[Theorem 3.1]{CNS}. In this case $A'\neq U$ since $\Sigma_2(A')=\Sigma_d(K)\neq S^3$.
\item The case when $K_q$ is knotted and $\Sigma_d(K_q)$ does not admit $SU(2)$ representations follows from  \cite[Theorem 4.1]{smalldehn}. To rule out the possibility that $A'=U$, we proceed by contradiction. Notice that the 2-fold cover of $\Sigma_d(K_q)$ branched over an unknot is homeomorphic to a sum $\Sigma_d(K_q)\# \Sigma_d(K_q)$. However, $\Sigma_d(K)$ is irreducible since $K$ is a prime knot (see \cite{plotnick}), which leads to a contradiction.
\end{enumerate}

If $\gcd(w,d)\neq 1$, set $d'=\tfrac{d}{\gcd(w,d)}$ and notice that $\gcd(w,d')=1$. The covering map $p:\Sigma_d(K)\to S^3$ then factors as $\Sigma_d(K)\to \Sigma_{d'}(K)\to S^3$, thus realizing $\Sigma_d(K)$ as a $\gcd(w,d)$-fold branched cover of $\Sigma_{d'}(K)$. By \Cref{surjections}, $\Sigma_{d'}(K)$ is an integer homology sphere and since $\gcd(w,d')=1$, the previous case establishes the existence of an irreducible representation $\rho:\pi_1\left(\Sigma_{d'}(K)\right)\to SU(2)$. The surjection from \Cref{surjections} gives the desired irreducible $SU(2)$ representation for $\Sigma_d(K)$.
\end{proof}

Knots with bridge number 2 and pretzel knots with odd parameters are examples of 2-periodic knots with unknotted quotient. In particular, every 2-bridge knot is $SU(2)$-simple (that is, traceless representations are binary dihedral and thus simple). Thus, our result is in stark contrast with a result of Sivek-Zentner \cite{menagerie} showing that 2-fold covers of $SU(2)$-simple knots of bridge number at most 3 do not admit any irreducible representations. \\

\begin{figure}[h]
\centering
\includegraphics[height=2.5cm]{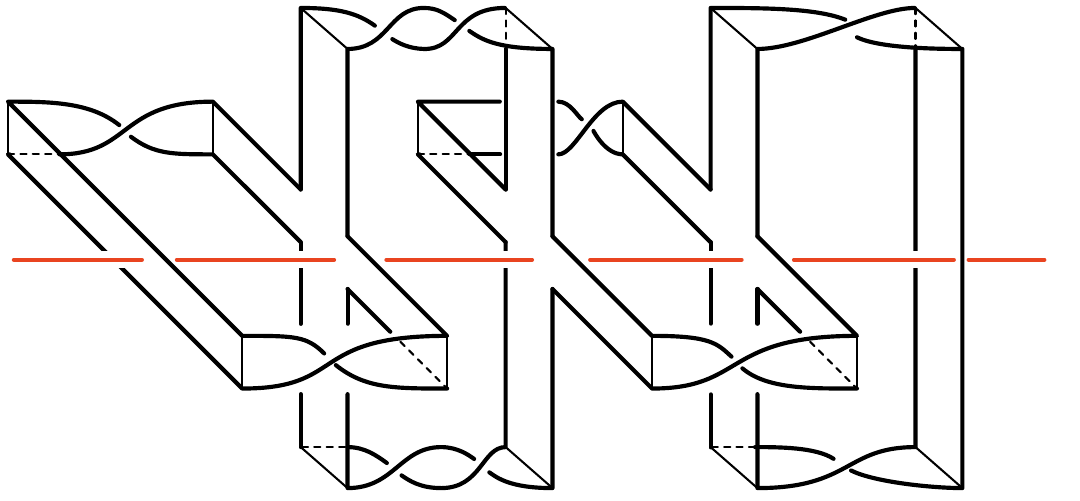} \qquad \includegraphics[height=2.5cm]{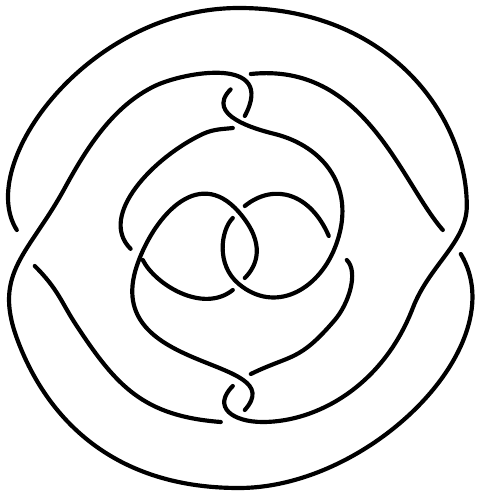} \qquad \includegraphics[height=2.5cm]{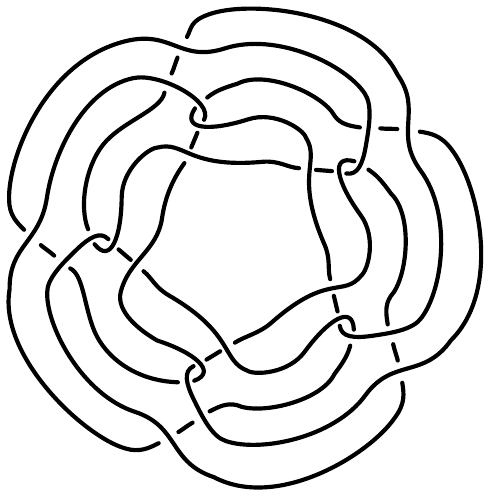}
\caption{The knot $K$ with even continued fraction expansion $23/16=[-2,2,-4,2]$. Left: A depiction of the 2-bridge knot and the periodicity axis. Center: A projection of the symmetric diagram of the knot $K$. Right: The knot $J$ from \Cref{2per-to-2fold} satisfying $\Sigma_5(K)=\Sigma_2(J)$.}\label{2-bridge}
\end{figure}

\begin{figure}[h]
\centering
\includegraphics[height=3cm]{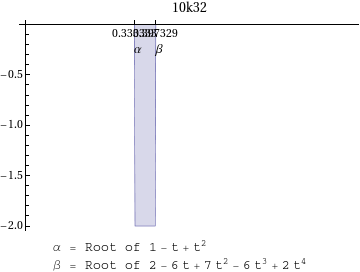}\qquad 
\includegraphics[height=3cm]{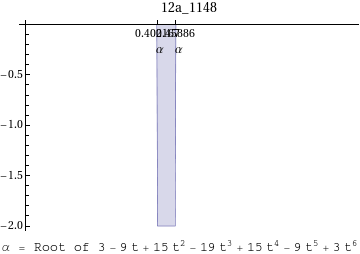}
\caption{The signature functions for $10_{32}$ (left) and $12a_{1148}$ (right). The images are taken directly from Knotinfo \cite{knotinfo}.}\label{thin-signature}
\end{figure}

To illustrate \Cref{2per-to-2fold} and its corollary, consider the 2-bridge knot $8_6$. It corresponds to the fraction $23/7$, has Alexander polynomial $\Delta_K(t)=2-6t+7t^2-6t^3+2t^4$, and has a 5-fold cover that is an integer homology sphere. To get a 2-periodic diagram of $8_6$, one uses the even continued fraction expansion $23/16=[-2,2,-4,2]$ and the `skewer' diagram shown on the left of \Cref{2-bridge}. The knot $J$ whose 2-fold cover is the 5-fold cover of $K$ is shown on the right of \Cref{2-bridge}. We remark that while \Crefparts{mainthm}{2per} shows that $\Sigma_5(K)$ admits an irreducible representation, the same result can be obtained by noticing that the signature function of $K$ is non-zero at a 5-fold root of unity.  Nevertheless, the signature criterion does not apply to the 2-bridge knots $10_{32}=69/19=[2,-2,2,2,-4]$ and $12a_{1148}=73/23$ since the support of the signature function is `thin' as shown in \Cref{thin-signature}. \\

To finish this section, we remark that the commuting trick used in the above proof of \Crefparts{mainthm}{2per} does not apply to invertible knots. This is because covers of the quotient theta graph commute only if all the degrees are 2. In the future, we hope to explore the remaining forms of symmetries: invertibility and both forms of amphichirality. The negative amphichiral case is particularly interesting, since these knots have finite order in algebraic concordance and therefore the signature function vanishes.

\section{The $\nu^{\sharp}$ invariant and $SU(2)$ representations}\label{instantons}
Because $SU(2)$-representation varieties of three manifolds are the set of critical points of Chern-Simons functionals, and these in turn are the generators of instanton Floer homology, the existence of $SU(2)$-representations can be obtained in terms of instanton homology groups or their chain complexes. In this section, we provide a criterion for the existence of a generator of the instanton chain complex that is not homologous to the trivial connection. We will use both the framed and the equivariant singular instanton groups. On one hand, the framed version can be defined for arbitrary $3$-manifold and many of its properties have been studied extensively, on the other hand, the equivariant theory has a unique property that enables us to conclude the existence of irreducible $SU(2)$-representation. Our proof thus combines the two theories.

For a knot $K$ in $S^3$, $\nu^{\sharp}(K)$ is a concordance invariant defined by Baldwin and Sivek \cite{baldwin2020concordance} using framed instanton homology with complex coefficients and based on appropriately chosen cobordism maps 
\[I^{\#}(X_n(K),{\nu_n}; \C): I^{\#}(S^3; \mathbb{C}) \to I^{\#}(S^3_n(K); \mathbb{C}),\] 
where the 4-manifold $X_n(K)$ is the trace of $n$-surgery.  More precisely, if $N(K)$ denotes the quantity
\begin{equation}\label{def:N}
N(K)=\min\{N\geq 0 \mid I^{\#}(X_n(K),{\nu_n}; \C)=0 \text{ for all } n\geq N\},\footnote{see the proof of \cite[Prop. 3.2]{baldwin2020concordance}}
\end{equation}
and $\overline{K}$ denotes the mirror image of $K$, then $\nu^{\sharp}(K)$ is the integer defined by 
\begin{equation}\label{def:nu}
\nu^{\sharp}(K) = N(K) - N(\overline{K}).
\end{equation}
\Cref{thm:nu} involving the $\nu^{\sharp}$ invariant is the main technical input we need from instanton Floer homology to prove the existence of irreducible $\SU(2)$ representations. Its proof relies on a careful analysis of the integers $n$ for which the $n$-trace cobordism maps vanish. However, to relate the cobordism map for framed instanton Floer homology with one for equivariant singular instanton Floer homology, we need to study a map in some cases different from Baldwin and Sivek's. Namely, we always fix the bundle to be trivial, but in \cite{baldwin2020concordance}, a possibly non-trivial bundle data $\nu_n$ is adopted. We remark that for $n\neq -1$, although $I^{\#}(X_n(K),{\nu_n}; \C)$ and $I^{\#}(X_n(K),\emptyset; \C)$ are not necessarily the same map, one is trivial exactly when the other one is. As we only work with complex coefficients in this paper, we will omit the coefficients from our notation from now on.

The following result has implicitly been proved in \cite{baldwin2020concordance}. For the self-containedness of the current paper, we present a proof here.
\begin{lemma}\label{trace_trivial}
	Let $K\subseteq S^3$ be a knot. For any integer $n$, the $n$-trace cobordism map $I^{\#}(X_n(K), \emptyset)$ is trivial if and only if one of the following happens.
\begin{itemize}
\item [(1)] $\nu^{\sharp}(K)\neq 0$\footnote{By \cite{baldwin2020concordance}, this assumption implies that $K$ is a V-shaped knot.} and $n\geq \nu^{\sharp}(K)$.
\item [(2)] $\nu^{\sharp}(K) = 0$, $K$ is V-shaped, and $n\geq -1$.
\item [(3)]\label{3} $\nu^{\sharp}(K) = 0$, $K$ is W-shaped, and $n\geq 1$.
\end{itemize}
\end{lemma}
\begin{proof} Let $F_n = I^{\#}(X_n(K), \emptyset; \C)$. The following is the exact sequence from \cite[Theorem 2.1 and Section 7.5]{scaduto} for the triple $(Y,Y_0,Y_1)=(S^3_{n+1}(K),S^3,S^3_{n}(K))$:
\begin{equation}\label{exact_triangle}
\xymatrix@C=3ex{
\cdots\ar[r]&I^\#(S^3_{n+1}(K),\mu)\ar[r]&I^\#(S^3)\ar[r]^{F_n\quad}&I^\#(S^3_{n}(K))\ar[r]&I^\#(S^3_{n+1}(K),\mu)\ar[r]& \cdots\\
}\end{equation}

To establish the desired triviality of $F_n$, notice that since $\dim_{\mathbb{C}}I^{\#}(S^3) = 1$, one can verify explicitly that
\begin{equation}\label{Fn=0}
F_n = 0~\Leftrightarrow~ \dim_{\mathbb{C}}I^{\#}(S^3_{n+1}(K);\mu) = \dim_{\mathbb{C}}I^{\#}(S^3_{n}(K)) + 1.
\end{equation}
For any knot $K\subseteq S^3$, and $n\neq 0$ we have 
\begin{align*}
\dim_{\mathbb{C}}I^{\#}(S^3_{n}(K))&=\dim_{\mathbb{C}}I^{\#}(S^3_{n}(K);\mu) &&\quad\text{(\cite[Theorem 1.12]{baldwin2020concordance2})}\\
\dim_{\mathbb{C}}I^{\#}(S^3_{n}(K))&=r_0(K)+|n-\nu^\sharp(K)| &&\quad\text{(\cite[Theorem 3.7]{baldwin2020concordance})}\end{align*}
and thus, as long as $n\neq 0,-1$, then $$ n\geq \nu^\sharp(K) ~\Rightarrow~ F_n = 0\qquad\text{and}\qquad n< \nu^\sharp(K) ~\Rightarrow~ F_n \neq 0.$$

In the case when $\nu^{\sharp}(K)\neq 0$ [Case (1)], \cite[Theorem 1.13]{baldwin2020concordance2} gives $\dim_{\mathbb{C}}I^{\#}(S^3_{0}(K),\mu) = \dim_{\mathbb{C}}I^{\#}(S^3_{0}(K)),$ and \cite[Theorem 3.7]{baldwin2020concordance} that $\dim_{\mathbb{C}}I^{\#}(S^3_{0}(K))=r_0(K)+|\nu^\sharp(K)|,$ so that the equivalences
	\[
		F_n = 0~\Leftrightarrow~ \dim_{\mathbb{C}}I^{\#}(S^3_{n+1}(K)) = \dim_{\mathbb{C}}I^{\#}(S^3_{n}(K)) + 1 \Leftrightarrow n\geq \nu^{\sharp}(K)
	\]
hold even for $n=0,-1$.

When $\nu^{\sharp}(K) = 0$, it remains to determine if the right hand side of \eqref{Fn=0} holds for $n=0,-1$. If $K$ is V-shaped [Case (2)], $r_0(K)= \dim_{\mathbb{C}}I^{\#}(S^3_{0}(K),\mu)=\dim_{\mathbb{C}}I^{\#}(S^3_{0}(K))+2$ (\cite[Definition 3.6]{baldwin2020concordance}, \cite[Theorem 1.13]{baldwin2020concordance2}). For the exceptional value $n=0$, we have  
\[
		\begin{aligned}
			\dim_{\mathbb{C}}I^{\#}(S^3_{1}(K),\mu)&=\dim_{\mathbb{C}}I^{\#}(S^3_{1}(K)) &&\quad\text{(\cite[Theorem 1.12]{baldwin2020concordance2})}\\
			&= \dim_{\mathbb{C}}I^{\#}(S^3_{0}(K)) + 1 &&\quad\text{(\cite[Theorem 3.7]{baldwin2020concordance})}
		\end{aligned}
	\]	
so that $F_0=0$. Now for the only exceptional value $n=-1$, we have
	\[
		\begin{aligned}
			\dim_{\mathbb{C}}I^{\#}(S^3_{0}(K),\mu) &= \dim_{\mathbb{C}}I^{\#}(S^3_{0}(K)) + 2 &&\quad\text{(\cite[Theorem 1.13]{baldwin2020concordance2})}\\
			&= \dim_{\mathbb{C}}I^{\#}(S^3_{-1}(K)) + 1 &&\quad\text{(\cite[Theorem 3.7]{baldwin2020concordance})}
		\end{aligned}
	\]
so that $F_{-1}=0$. We conclude that, in this case, $F_n=0$ if and only if $n\geq -1$.
	
	Lastly, the case when $K$ is W-shaped [Case (3)]  can be argued similarly. In this case $r_0(K)= \dim_{\mathbb{C}}I^{\#}(S^3_{0}(K),\mu)=\dim_{\mathbb{C}}I^{\#}(S^3_{0}(K))-2$ (\cite[Definition 3.6]{baldwin2020concordance}, \cite[Theorem 1.13]{baldwin2020concordance2}). For the exceptional value $n=0$, we have  
\[
		\begin{aligned}
			\dim_{\mathbb{C}}I^{\#}(S^3_{1}(K),\mu)&=\dim_{\mathbb{C}}I^{\#}(S^3_{1}(K)) &&\quad\text{(\cite[Theorem 1.12]{baldwin2020concordance2})}\\
			&= \dim_{\mathbb{C}}I^{\#}(S^3_{0}(K),\mu) + 1 &&\quad\text{(\cite[Theorem 3.7]{baldwin2020concordance})}\\
			&= \dim_{\mathbb{C}}I^{\#}(S^3_{0}(K)) - 1 &&\quad\text{(\cite[Theorem 1.13]{baldwin2020concordance2})}
		\end{aligned}
	\]	
and for the exceptional value $n=-1$, we have
	\[
		\begin{aligned}
			\dim_{\mathbb{C}}I^{\#}(S^3_{0}(K),\mu) &= \dim_{\mathbb{C}}I^{\#}(S^3_{0}(K)) - 2 &&\quad\text{(\cite[Theorem 1.13]{baldwin2020concordance2})}\\
			&= \dim_{\mathbb{C}}I^{\#}(S^3_{-1}(K)) - 1 &&\quad\text{(\cite[Theorem 3.7]{baldwin2020concordance})}
		\end{aligned}
	\]
	Thus, we conclude that, in this case, $F_n=0$ if and only if $n\geq 1$.

\end{proof}

We now introduce the next instanton homology: the \emph{equivariant singular instanton homology} as developed by Daemi and Scaduto in \cite{DS-singular, DS-clasp}. For a pair $(Y,K)$ with $Y$ an oriented integer homology sphere and $K\subset Y$ a knot, the authors consider $SU(2)$-connections on $Y$ with singularities of order 4 along the knot $K$. They define a chain complex $\left(C_*(Y, K),d\right)$ over $\C$ with generators given, roughly, by gauge equivalence classes of $SU(2)$-connections which, in addition to being singular, are also flat and irreducible (see \cite[Section 2.2]{DS-singular}). The differential counts singular instantons on $[0,1]\times Y$ that restrict to prescribed flat connections on each end (\cite[Sections 2.5 and 3.1]{DS-singular}). The underlying $\C$-vector spaces of the equivariant singular instanton complex $\left(\widetilde{C}(Y, K),\tilde{d}\,\right)$ take the form
\begin{equation}\label{eq: equivariant singular chain complex}
\widetilde{C}_*(Y, K)= C_* (Y, K)\oplus C_{*-1}(Y, K) \oplus \C[\theta],\end{equation}
where $C_{*-1}(Y, K)$ is isomorphic as a $\C$-vector space to $C_* (Y, K)$ but differs by a grading-shift, and $\C[\theta]$ is generated by the unique reducible singular $\SU(2)$ connection $\theta$. The differential $\tilde{d}$ takes the form \begin{equation}
    \tilde{d}= \begin{pmatrix}
d & 0 & 0 \\
v & -d & \delta_2 \\  
\delta_1 & 0 & 0 
\end{pmatrix}.\end{equation}
Here, $\delta_1,\delta_2$ are chain maps that count instantons in such a way as to incorporate the reducible connection $\theta$, and $v$ is an anti chain homotopy between the composition $\delta_2\circ \delta_1$ and the zero map (cf. \cite[Proposition 3.16]{DS-singular}). See \cite[Section 3.2]{DS-singular} for the definition of the maps $\delta_1,\delta_2$, \cite[Section 3.3.2]{DS-singular} for the definition of $v$, and \cite[Section 3.4]{DS-singular} for a complete overview of the complex $(\widetilde{C}(Y, K),\tilde{d}\,)$ and its homology groups $\widetilde{I}(Y, K)$.

Now, assume $Y$ is a closed oriented $3$-manifold obtained from $S^3$ via integral surgeries along a link $L$ with respect to a linking matrix $A$. This description of $Y$ induces a $2$-handle cobordism $W:S^3\to Y$. For $W$ to induce a map on \emph{equivariant singular instanton homology}, we first have to construct a surface as follows: let $U\subseteq S^3$ be an unknotted curve that is split from $L$, and form the annulus $S=[0,1]\times U$. This construction gives a cobordism of pairs $(W,S): (S^3,U)\to (Y,U)$. If the matrix $A$ is unimodular and negative definite, then $(W,S)$ is a negative definite pair in the sense of \cite[Definition 2.33]{DS-singular}, and induces well-defined maps $\widetilde{C}(W, S) \colon \widetilde{C}(S^3, U) \to \widetilde{C}(Y, U)$ and $\tilde{I}(W, S) \colon \tilde{I}(S^3, U) \to \tilde{I}(Y, U)$. The following proposition relates the maps induced by cobordisms $W:S^3\to Y$ to both the \emph{framed} and the \emph{equivariant singular} instanton theories. 

\begin{proposition}\label{sharp-tilde} The cobordism map $I^{\#}(W,\emptyset)$ in framed instanton Floer homology is trivial if and only if the cobordism map $\widetilde{I}(W, S)$ in equivariant singular instanton Floer homology is trivial.
\end{proposition}

\begin{proof}
We have the following commutative diagram.

\begin{equation}\label{diag:1}
\xymatrix{
\widetilde{I}(S^3, U)\ar[d]^{\widetilde{I}(W, S)}\ar[r]^\iso & I^{\natural} (S^3, U)\ar[d]^{I^\natural(W, S)}\ar[r]^\iso & \KHI(S^3, U)\ar[d]^{{KHI(W, S)}} \ar[r]^{\iso} & I^{\#}(S^3) \ar[d]^{I^{\#}(W)} \\
\widetilde{I}(Y, U) \ar[r]^{\iso} & I^{\natural} (Y, U) \ar[r]^{\iso} & \KHI(Y, U) \ar[r]^{\iso} & I^{\#}(Y)
}
\end{equation}

The commutativity of the left square of Diagram~\ref{diag:1} comes from \cite[Theorem~8.9]{DS-singular} while the commutativity of the middle square comes from the proof of \cite[Corollary 4.5]{DLVW-ribbon}.

For the right square, following \cite[Definition~7.13]{kronheimer2010knots}, we first note that for any $3$-manifold $Z$, $\KHI(Z,K) = I_*(Y_1|\bar R)_w$, where $Y_1$ is a specific closure of the sutured manifold $(Z \setminus N(K), \Gamma_\mu)$, where $\Gamma_\mu$ is the meridional suture. The closure $Y_1$ is formed by gluing $F\times S^1$ to the knot complement $Z \setminus N(K)$ so that ${p} \times S^1$ is attached to the meridian of $K$ and $\partial F \times {q}$ is glued to any chosen longitude of $K$, where $F$ is a genus-1 surface with one boundary component. Since $\bar R$ has genus 1, $I_*(Y_1|\bar R)_w$ can be characterized as the $+2$ eigenspace of $\mu(y)$-action on $I_*(Y)$. Following \cite[Definition~2.2, Remark~2.3]{lidman2022framed}, this $+2$ eigenspace of $\mu(y)$ is then identified with $I^{\#}(Z)$. This proves that for any $3$-manifold $Z$, $\KHI(Z, U) \iso I^{\#}(Z)$. To prove the claim that this isomorphism is natural with respect to the cobordism map, it suffices to observe that since the associated closure for both manifolds is $S^3 \# T^3$ and the cobordism is given by 2-handle attachment along a link $L \subset S^3 \# T^3$, the associated 4-manifold in which we count finite energy instantons is the same. Hence, the cobordism maps in these two categories naturally agree. This completes the proof of the claim about the commutativity of \eqref{diag:1}.
\end{proof}

\begin{theorem}\label{thm:nu} Let $A$ be an $n\times n$ unimodular matrix that is symmetric and negative definite. Suppose $L = L_1\cup\cdots L_n$ is a framed link with linking matrix $A$. Suppose there exists an $i\in\{1,2,\cdots, n\}$ such that
\begin{enumerate}
\item $\nu^\sharp(L_i)\neq 0$ and $a_{ii}\geq \nu^\sharp(L_i)$, or
\item $\nu^\sharp(L_i)=0$, $L_i$ is $V$-shaped, and $a_{ii}\geq -1$, or
\item $\nu^\sharp(L_i)=0$, $L_i$ is $W$-shaped and $a_{ii}\geq 1$.
\end{enumerate}Then, the 3-manifold $Y$ with surgery description given by the framed link $L$ is an integer homology sphere whose fundamental group admits an irreducible $\SU(2)$-representation.
\end{theorem}

\begin{proof}
Let $W$ be the cobordism from $S^3$ to $Y$ corresponding to the above surgery description. It is well known that the cobordism map $I^\# (W) \colon I^\# (S^3) \to I^\#(Y)$ is invariant under reordering of handle attachments. Hence, reordering the handles if necessary, we can assume that $a_{11}\geq\nu^\sharp(L_1)$, and decompose $W$ as 
\[
W = X_{a_{11}(L_1)}\cup W',
\]
where $X_{a_{11}(L_1)}$ is the trace cobordism from $S^3$ to $S^3_{a_{11}}(L_1)$, and $W'$ is the result of attaching the remaining 2-handles to the outgoing boundary component of $X_{a_{11}(L_1)}$. The naturality and functoriality of the instanton cobordism map then implies
\[
I^{\#}(W) = I^{\#}(W')\circ I^{\#}(X_{a_{ii}(L_i)}).
\]
The hypotheses on $a_{11}$ together with \Cref{trace_trivial} imply that $I^\#(X_{a_{11}(L_1)})=0$, and so by the composition law $I^{\#}(W)=0$. By \Cref{sharp-tilde}, the cobordism map $\widetilde{I}(W; S)$ is also trivial:
\begin{equation}\label{eq: I-tilde is trivial}
	\tilde{I}(W,S)= 0.
\end{equation}
Our next goal is to relate this triviality to the existence of irreducible representations. 

Following \cite[Section 3]{DS-singular}, the negative definite pair $(W,S)$ induces a chain map
\begin{equation}\label{eq: cobordism map in equivariant singular instanton Floer homology}
	\widetilde{CI}(W,S) = \left(
	\begin{matrix}
		\lambda & 0 & 0\\
		\mu & \lambda & \Delta_2\\
		\Delta_1 & 0 & z
	\end{matrix}
\right): \widetilde{C}(S^3,U)\to \widetilde{C}(Y,U)
\end{equation}
where $\lambda$, $\mu$, $\Delta_1$, and $\Delta_2$ can be regarded as 4-dimensional analogues of $d$, $v$, $\delta_1$, and $\delta_2$. We remark that $z$, the coefficient from the canonical reducible generator $\theta_{(S^3,U)}$ to $\theta_{(Y,U)}$ is a non-zero complex number.

We now prove that $\pi_1(Y)$ must admit an irreducible $SU(2)$ representation. Assume the contrary, i.e., all $SU(2)$-representations of $\pi_1(Y)$ are reducible. We know that reducible $SU(2)$-representations are abelian, and since $Y$ is an integral homology sphere, there is a unique (trivial) representation $\rho: \pi_1(Y)\to SU(2)$. We can then write $Y \setminus U = Y \# (S^3 \setminus U)$ and conclude that there is a unique $SU(2)$-representation $\theta_Y$ of $\pi_1(Y\backslash U)$ that sends the meridian of $U$ to $i\in SU(2)$. Note this unique representation must be reducible and by \cite[Proposition 2.3]{DS-singular} it is isolated and non-degenerate. Thus we can build the associated equivariant chain complex with no irreducible generators, and thus
\[
\widetilde{C}(Y,U) = 0\oplus 0\oplus <\theta_Y>.
\]
Similar argument applies to $S^3$ as well, and we have
\[
\widetilde{C}(S^3,U) = 0\oplus 0\oplus <\theta_{S^3}>.
\]
Hence, the cobordism chain map will be reduced to 
\[
\widetilde{CI}(W,S) (\theta_{S^3}) = z\cdot \theta_Y\text{ where }z\text{ is non-zero}.
\]
This leads to the fact that $\tilde{I}(W,S)\neq 0$, which directly contradicts (\ref{eq: I-tilde is trivial}).

\end{proof}

\begin{remark} In the ideal case when perturbations are not necessary to define the complex $C_*(Y,U)$, the triviality of the cobordism map $\tilde{I}(W,S)$ would imply the existence of $\alpha$ with $d(\alpha)=0$ with $\delta_1(\alpha) \neq 0$. Following \cite[Proposition 4.15]{DS-singular}, this would imply $h(Y, U)>0$. We expect that the nonvanishing of $h(Y, U)$ implies the nonvanishing of the Fr\o yshov invariant $h(Y)$ \cite[Definition 3]{Froyshov-h}.
\end{remark}


\section{A construction of branched covers admitting $SU(2)$ representations}\label{sicup}
In this section we use the instanton setup from \Cref{instantons}, and show that covers of $S^3$ branched along closures of some special tangles admit a irreducible $SU(2)$-representations. We first define SICUP matrices and their properties, then review equivariant surgery descriptions of branched covers, and finally prove \Crefparts{mainthm}{nu}. 

\subsection{A quick detour into matrices} This subsection is an attempt to elucidate some basic properties of the types of matrices considered in \Crefparts{mainthm}{nu}, and in the construction of families of knots whose $d$-fold covers admit $SU(2)$-representations. In particular, we classify all SICUP matrices of size at most 5.

\begin{definition}\label{def:circulant} A circulant matrix is a $d\times d$ square matrix $C$ in which all rows are composed of the same elements and each row is rotated one element to the right relative to the preceding row. 
$\displaystyle 
C=\begin{bmatrix}
c_{1}&c_{2}&\cdots &c_{d-1}&c_{d}\\
c_{d}&c_{1}&c_{2}&&c_{d-1}\\
\vdots &c_{d}&c_{1}&\ddots &\vdots \\
c_{3}&&\ddots &\ddots &c_{2}\\
c_{2}&c_{3}&\cdots &c_{d}&c_{1}\\
\end{bmatrix}$

\end{definition}

As mentioned at the end of the proof of \Cref{surgery_cover}, our linking matrices are not just circulant, but they are symmetric and their entries are integers as well. The additional criteria below are requirements from \Cref{thm:nu}.

\bdefn\label{defn: sicup matrices}
A $d\times d$ matrix $A$ is called an {\bf SICUP matrix} if the following hold.
\begin{itemize}
	\item ({\bf S}ymmetric) $A=A^T$
	\item ({\bf I}ntegral) $A\in M_{d\times d}(\Z)$, that is, $A$ has integral entries.
	\item ({\bf C}irculant) For any pair $(i,j)\in\{1,\cdots,d\}^2$, we have $a_{i,j} = a_{i+1,j+1}$ (taking indices $\mod d$).
	\item ({\bf U}nimodular) We have ${\rm det}(A) = 1$.
	\item ({\bf P}ositive definite) The matrix $A$ is positive definite.
\end{itemize}
Let $\mathcal{M}(d)$ be the set of all SICUP matrices of size $d$.
\edefn

The size of the matrix will be related to the degree of a cover of $S^3$. Since a $2r$-cover is the 2-fold cover of an $r$-fold cover, the existence of $SU(2)$-representations is guaranteed by \Cref{2fold-thm}. We thus restrict our attention to the case of odd size $d=2r+1$ for $r>0$ and include some basic facts regarding SICUP matrices. We omit the proofs but point the reader to \cite[Chapter 3]{gray} for a nice review.

\begin{remark}\label{prop:circulant} Let $C$ be a $d\times d$ SICUP matrix with $d=2r+1$, $r\in\N$. Then we have:
\begin{enumerate}
\item If $\vec{u}$ denotes the vector of ones, then $\vec{u}$ is an eigenvector for $C$ with eigenvalue $\lambda_1=c_1+2(c_2+\ldots c_r)$. Positivity and unimodularity imply that $\lambda_1=1$.
\item There is an eigenvalue $\lambda_j$ of multiplicity 2 for $j=2,\ldots,r$. This allows us to write $\det(C)=\lambda_1(\lambda_2\cdot\ldots\cdot\lambda_r)^2=1$.
\end{enumerate}
\end{remark}

Simple computations show that for $d=1,2,3,4$, the set $\mathcal{M}(d)$ has the $d\times d$ identity matrix as its only element. For $d=5$ we have the following intriguing relationship between SICUP matrices and Pell's equation.  

\begin{proposition}\label{prop: SICUP and Pell's eq}
There is a one-to-one correspondence between the set of $5\times 5$ integral symmetric matrices $\mathcal{M}(5)$ and the set of solutions $(a,b)$ to the generalized Pell's equation $a^{2}-5b^{2}=4$ such that $a\equiv 2~(\text{mod }5)$.
\end{proposition} 
\begin{proof} Notice that any $5\times 5$ integral symmetric circulant matrix can be parametrized by three integers $x$, $l$, and $m$:
\[
\renewcommand{\arraystretch}{1.25}
A(x,l,m) = \left[
\begin{matrix}
	x & l & m & m & l\\
	l & x & l & m & m\\
	m & l & x & l & m\\
	m & m & l & x & l\\
	l & m & m & l & x\\
\end{matrix}
\right].\]
A straightforward computation shows that the eigenvalues of $A$ are given by
\begin{align*}
\lambda_1&=x+2l+2m\\
\lambda_2&=\tfrac{1}{2}\left((2x-l-m)+\sqrt{5}(l-m)\right)\\
\lambda_3&=\tfrac{1}{2}\left((2x-l-m)+\sqrt{5}(m-l)\right),
\text{ so that}\\
\det(A)&=\lambda_1\left(\lambda_2\lambda_3\right)^2, \text{ with}\\
\lambda_2\lambda_3&=\tfrac{1}{2}\left((2x-l-m)^2- 5(l-m)^2\right).
\end{align*}
By \Cref{prop:circulant}, ${\rm det}(A)=1$ and $\lambda_1=1$, forcing $\lambda_2\lambda_3= 1$. Thus, for $A$ to be an element of $\mathcal{M}(5)$ we must have 
\begin{equation*}
(2x-l-m)^2- 5(l-m)^2=4.
\end{equation*}
This immediately defines a map
\begin{align*}
\Phi: \mathcal{M}(5)&\longrightarrow\Pell(5,4)=\{(a,b)\in\Z^2 \mid\ a^{2}-5b^{2}=4\}\\
A(x,l,m)&\longrightarrow(2x-l-m,l-m).
\end{align*}
To show that this map is a bijection, first notice that if $(a,b)\in\Pell(5,4)$, then $a\equiv 2 \pmod{5}$. An explicit inverse for $\Phi$ is obtained as 
\begin{align*}
\Phi^{-1}:\Pell(5,4)&\longrightarrow\mathcal{M}(5)\\
(a,b)&\longrightarrow A\left(\frac{2a+1}{5},\frac{2-a+5b}{10},\frac{2-a-5b}{10}\right),
\end{align*}
where we have used the fact that $x=1-2l-2m$.
\end{proof}

\subsection{Equivariant surgery descriptions of branched covers}

There are many different ways of describing branched covers of $S^3$ as surgery on framed links. Among those are the algorithms of Akbulut-Kirby using Seifert surfaces \cite{akbulut-kirby}, or those whose input is a surgery diagram for the branching set as in \cite[Chapters 6C, 10C]{rolfsen} or \cite[Theorem 3]{goldsmith}. A careful analysis of the latter algorithm shows that surgery descriptions for all the covers can be obtained from the same set of unknotted circles. That is, there is no restriction on the degree of the cover. Since we don't need the result in its greatest generality, we content ourselves with a description for a specific choice of covering degree $d$. Equivariant surgery descriptions of branched covers establish the technical framework for our SICUP matrix approach, and the constructions that will appear in the next subsection.

\begin{definition} Let $\boldsymbol\gamma=\gamma_1\sqcup\ldots\sqcup\gamma_u\subset S^3\setminus N(K)$ be a $u$-component unlink. If $K'$ is the knot obtained from $K$ after performing $\pm 1$ surgery on $S^3$ along $\boldsymbol\gamma$, then $K'$ is said to be obtained from $K$ by a sequence of $u$ twists.

If $K'$ is the unknot, we say that the link $\boldsymbol\gamma$ untwists $K$. If in addition $lk(K,\gamma_i)=0$ for $i=1,\ldots,u$, we say $\boldsymbol\gamma$ unknots $K$.
\end{definition}

We claim no originality in the following statement; we simply include it to make the constructions of the next section clearer.

\begin{proposition}\label{surgery_cover} Let $d\geq 2$ be an integer, and $\boldsymbol\gamma=\gamma_1\sqcup\ldots\sqcup\gamma_u$ a link that untwists $K$ and satisfies $lk(K,\gamma_i)\equiv 0 \pmod{d}$. Then $\Sigma=\Sigma_d(K)$ can be obtained from $S^3$ by equivariant surgery on the link $L$ obtained as the preimage of $\boldsymbol\gamma$ in $S^3=\Sigma_d(U)$. Moreover, the link $L$ satisfies the following:
\begin{enumerate}
\item The link $L$ is invariant under the set of covering transformations of $\Sigma$.
\item The framing curves that determine the surgery are also invariant under the covering transformations.
\item The linking matrix of $L$ is symmetric and has circulant blocks.
\end{enumerate}
\end{proposition}

\begin{proof} To prove the result, it will be enough to construct $\Sigma_d(K)$ via surgery on $\Sigma_d(U)=S^3$. If $p:\Sigma_d(U)\to S^3$ denotes the covering map, then $\Sigma_d(K)$ will be described as a filling of the covering space $p^{-1}\left(S^3\setminus N(U\sqcup \boldsymbol\gamma)\right)\to S^3\setminus N(U)$. The hypothesis $lk(K,\gamma_i)=lk(U,\gamma_i)\equiv 0 \pmod{d}$ implies that the preimage of $\gamma_i$ in $\Sigma_d(U)$ is a link with $d$-components. Let $f$ be a generator of the group of covering transformations, and denote by $\gamma^{(j)}_i$ the $j$-th lift of $\gamma_i$. Choose the numbering so that $f\left(\gamma^{(j)}_i\right)=\gamma^{(j+1)}_i$, and more generally, $f^j\left(\gamma^{(1)}_i\right)=\gamma^{(j+1)}_i$. Choose a meridian-longitude pair $(\mu_i,\lambda_i)$ for $\gamma_i$ so that any curve $\eta_i\in\partial N(\gamma_i)$ representing  $n_i\mu_i+\lambda_i$ is a framing curve for the surgery.  Denote by $\eta_i^{(j)}$ the lift of $\eta_i$ that lies in $\partial N\left(\gamma^{(j)}_i\right)$. Our conventions are such that, just as with the lifts of the $\gamma_i$'s, $f^j\left(\eta^{(1)}_i\right)=\eta^{(j+1)}_i$. From all of this, it follows that $\Sigma_d(K)$ is the filling of $\Sigma_d(U)\setminus p^{-1}(\boldsymbol\gamma)$ that caps off the curves $\eta_i^{(j)}$ with the meridional disks of the attached solid tori. Invariance under covering transformations thus follows.

As for the linking matrix, assume first that $u=1$. In this case, the linking matrix $\mathcal{L}$ is a $d\times d$ symmetric matrix with $(k,k)$ entry given by the framing of $\gamma^{(k)}$, that is,
\begin{equation}
l_{kk}=lk(\gamma^{(k)},\eta^{(k)})=lk(f^{k-1}\left(\gamma^{(1)}\right),f^{k-1}\left(\eta^{(1)}\right))=lk(\gamma^{(1)},\eta^{(1)}).
\end{equation}
Since the surgery description is invariant under covering transformations and $\gamma^{(k)}$ is the image of $\gamma^{(1)}$ under a covering transformation, all the diagonal entries of $\mathcal{L}$ have the same value $l_{11}$. Furthermore, the $(1,1)$-entry of $\mathcal{L}$ satisfies $l_{11}=n-\sum_{j=2}^d l_{1j}$ (see \cite[Lemma 2.3]{homo} or \cite[pg. 140]{faces}).
Since $\mathcal{L}$ is symmetric, the $(k,j)$-entry and the $(j,k)$-entry are equal, and so it is enough to compute $l_{kj}$ for $j>k$. In this case we have $$l_{kj}
=lk\left(\gamma^{(k)},\gamma^{(j)}\right)=lk\left(f^{k-1}\left(\gamma^{(1)}\right),f^{j-k}f^{k-1}\left(\gamma^{(1)}\right)\right)=lk\left(\gamma^{(1)},f^{j-k}\left(\gamma^{(1)}\right)\right)=l_{1,j-k+1}
.$$ 
In conclusion, the entries of the first row completely determine the entire matrix, where the $k$-th row is obtained from the first row by shifting its entries $k$ positions. Matrices with this property are called circulant (for a precise definition see \Cref{def:circulant} below).

If $u>1$, then the linking matrix is a $u\times u$ block matrix, with $d\times d$ blocks. Denote by $\mathcal{L}_{rs}$ the $(r,s)$-block of $\mathcal{L}$. The argument for the case $u=1$ applies (with appropriate modifications) to conclude that each block $\mathcal{L}_{rr}$ is symmetric circulant. For $s\neq r$ and $j>k$ we have:
$$\left(\mathcal{L}_{r,s }\right)_{kj}=lk\left(\gamma_r^{(k)},\gamma_s^{(j)}\right)=lk\left(\gamma_r^{(1)},f^{j-k}\left(\gamma_s^{(1)}\right)\right)=\left(\mathcal{L}_{r,s }\right)_{1, j-k+1}.$$
Thus, each block is circulant (but not necessarily symmetric circulant unless $r=s$).

\end{proof}


\begin{exmp}\label{exmp:10_132} For $K=10_{132}$, the cover $\Sigma=\Sigma_3(K)$ is an integer homology 3-sphere for which the existence of an irreducible representation $\rho:\pi_1(\Sigma)\to SU(2)$ is established by applying \Cref{thm:nu}. The equivariant surgery presentation $L=L_1\sqcup L_2\sqcup L_3$ from \Cref{10_132} has linking matrix $-I_{3\times 3}$, and each component $L_i$ of the surgery link is an unknot. Thus, to apply the theorem we need a different surgery presentation. The 2-component surgery description $L'_1\sqcup L'_2$ from \Cref{10_132-nu} is obtained after blowing down $L_3$. It has linking matrix $A=-I_{2\times 2}$, but each component $L'_i$ has the knot type of $5_2$\footnote{It is important to notice that by $5_2$ we mean the twist knot with a negative clasp}. Using the notation from \cite[Theorem 1.13]{baldwin2020concordance}, this is the knot $K_3$, and so $\nu^\#(L'_1)=-1=a_{11}$. 

\begin{figure}[h]
\centering
\def\svgwidth{0.25\textwidth}
\begingroup%
  \makeatletter%
  \providecommand\color[2][]{%
    \errmessage{(Inkscape) Color is used for the text in Inkscape, but the package 'color.sty' is not loaded}%
    \renewcommand\color[2][]{}%
  }%
  \providecommand\transparent[1]{%
    \errmessage{(Inkscape) Transparency is used (non-zero) for the text in Inkscape, but the package 'transparent.sty' is not loaded}%
    \renewcommand\transparent[1]{}%
  }%
  \providecommand\rotatebox[2]{#2}%
  \newcommand*\fsize{\dimexpr\f@size pt\relax}%
  \newcommand*\lineheight[1]{\fontsize{\fsize}{#1\fsize}\selectfont}%
  \ifx\svgwidth\undefined%
    \setlength{\unitlength}{342.99212598bp}%
    \ifx\svgscale\undefined%
      \relax%
    \else%
      \setlength{\unitlength}{\unitlength * \real{\svgscale}}%
    \fi%
  \else%
    \setlength{\unitlength}{\svgwidth}%
  \fi%
  \global\let\svgwidth\undefined%
  \global\let\svgscale\undefined%
  \makeatother%
  \begin{picture}(1,0.79338843)%
    \lineheight{1}%
    \setlength\tabcolsep{0pt}%
    \put(0,0){\includegraphics[width=\unitlength,page=1]{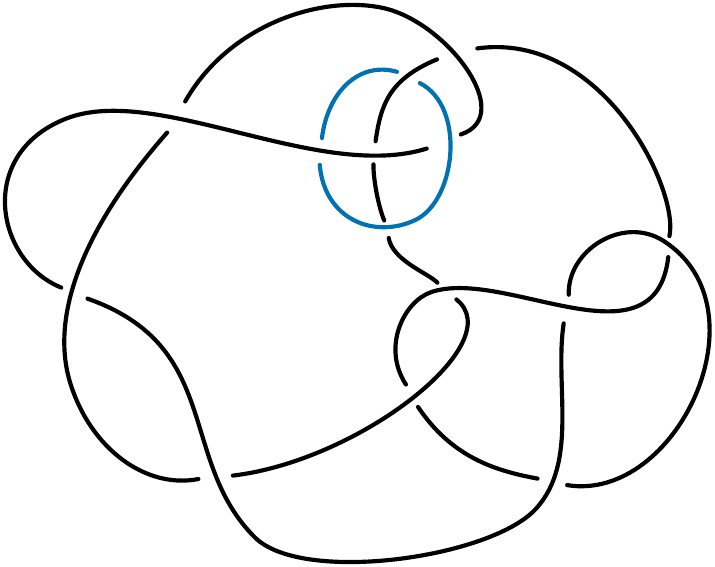}}%
    \put(0.1,0.1){\makebox(0,0)[rt]{\lineheight{1.25}\smash{\begin{tabular}[t]{r}$K$\end{tabular}}}}%
    \put(0.46114862,0.4){\color[rgb]{0,0.44705882,0.69803922}\makebox(0,0)[rt]{\lineheight{1.25}\smash{\begin{tabular}[t]{r}$-1$\end{tabular}}}}%
  \end{picture}%
\endgroup%
\qquad
\def\svgwidth{0.25\textwidth}
\begingroup%
  \makeatletter%
  \providecommand\color[2][]{%
    \errmessage{(Inkscape) Color is used for the text in Inkscape, but the package 'color.sty' is not loaded}%
    \renewcommand\color[2][]{}%
  }%
  \providecommand\transparent[1]{%
    \errmessage{(Inkscape) Transparency is used (non-zero) for the text in Inkscape, but the package 'transparent.sty' is not loaded}%
    \renewcommand\transparent[1]{}%
  }%
  \providecommand\rotatebox[2]{#2}%
  \newcommand*\fsize{\dimexpr\f@size pt\relax}%
  \newcommand*\lineheight[1]{\fontsize{\fsize}{#1\fsize}\selectfont}%
  \ifx\svgwidth\undefined%
    \setlength{\unitlength}{342.99212598bp}%
    \ifx\svgscale\undefined%
      \relax%
    \else%
      \setlength{\unitlength}{\unitlength * \real{\svgscale}}%
    \fi%
  \else%
    \setlength{\unitlength}{\svgwidth}%
  \fi%
  \global\let\svgwidth\undefined%
  \global\let\svgscale\undefined%
  \makeatother%
  \begin{picture}(1,0.79338843)%
    \lineheight{1}%
    \setlength\tabcolsep{0pt}%
    \put(0,0){\includegraphics[width=\unitlength,page=1]{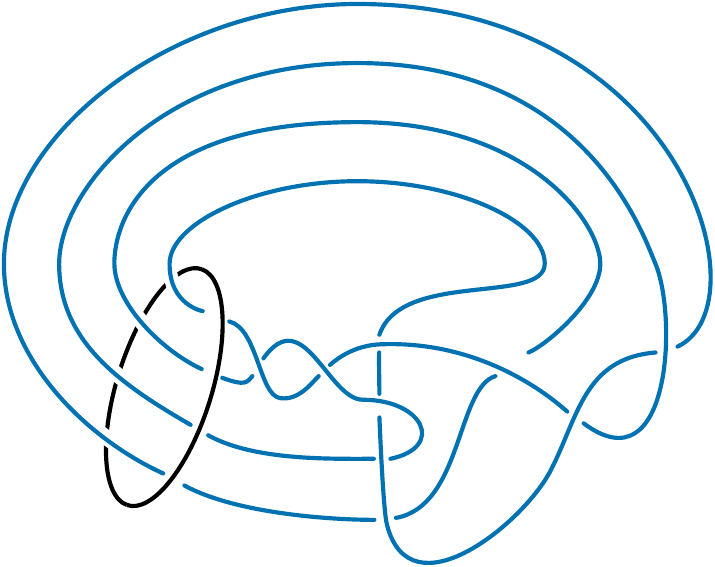}}%
    \put(0.14324557,0){\makebox(0,0)[rt]{\lineheight{1.25}\smash{\begin{tabular}[t]{r}$K$\end{tabular}}}}%
    \put(0.51818955,0.8){\color[rgb]{0,0.44705882,0.69803922}\makebox(0,0)[rt]{\lineheight{1.25}\smash{\begin{tabular}[t]{r}$-1$\end{tabular}}}}%
  \end{picture}%
\endgroup%
\qquad
\def\svgwidth{0.2125\textwidth}
\begingroup%
  \makeatletter%
  \providecommand\color[2][]{%
    \errmessage{(Inkscape) Color is used for the text in Inkscape, but the package 'color.sty' is not loaded}%
    \renewcommand\color[2][]{}%
  }%
  \providecommand\transparent[1]{%
    \errmessage{(Inkscape) Transparency is used (non-zero) for the text in Inkscape, but the package 'transparent.sty' is not loaded}%
    \renewcommand\transparent[1]{}%
  }%
  \providecommand\rotatebox[2]{#2}%
  \newcommand*\fsize{\dimexpr\f@size pt\relax}%
  \newcommand*\lineheight[1]{\fontsize{\fsize}{#1\fsize}\selectfont}%
  \ifx\svgwidth\undefined%
    \setlength{\unitlength}{342.99212598bp}%
    \ifx\svgscale\undefined%
      \relax%
    \else%
      \setlength{\unitlength}{\unitlength * \real{\svgscale}}%
    \fi%
  \else%
    \setlength{\unitlength}{\svgwidth}%
  \fi%
  \global\let\svgwidth\undefined%
  \global\let\svgscale\undefined%
  \makeatother%
  \begin{picture}(1,1)%
    \lineheight{1}%
    \setlength\tabcolsep{0pt}%
    \put(0,0){\includegraphics[width=\unitlength,page=1]{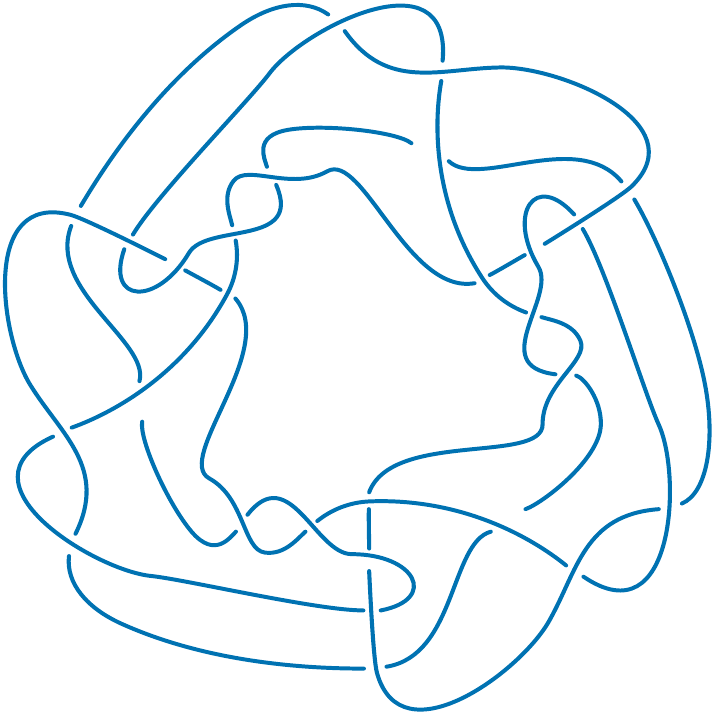}}%
    \put(1.1,0.2){\color[rgb]{0,0.44705882,0.69803922}\makebox(0,0)[rt]{\lineheight{1.25}\smash{\begin{tabular}[t]{r}$-1$\end{tabular}}}}%
   \put(0.45,1.01){\color[rgb]{0,0.44705882,0.69803922}\makebox(0,0)[rt]{\lineheight{1.25}\smash{\begin{tabular}[t]{r}$-1$\end{tabular}}}}%
    \put(0.1,0.1){\color[rgb]{0,0.44705882,0.69803922}\makebox(0,0)[rt]{\lineheight{1.25}\smash{\begin{tabular}[t]{r}$-1$\end{tabular}}}}%
  \end{picture}%
\endgroup%
\caption{The knot $K=10_{132}$ is hyperbolic and its 3-fold cover is an integer homology sphere.}\label{10_132}
\end{figure}
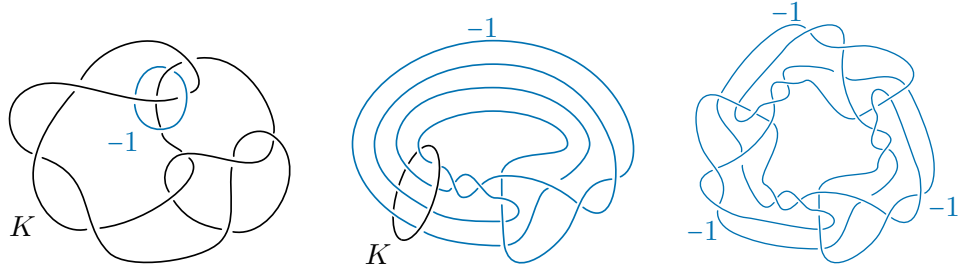

\begin{figure}[h]
\centering
\def\svgwidth{0.25\textwidth}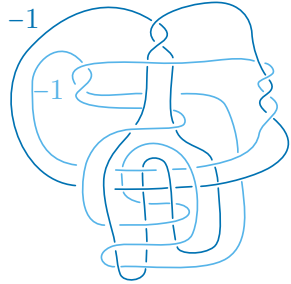
\caption{
A surgery description for $\Sigma_3(10_{132})$ that satisfies the hypotheses of \Cref{thm:nu}.}\label{10_132-nu}
\end{figure}
\end{exmp}

\subsection{Constructing branched covers admitting $SU(2)$ representations}

For a fixed $d$, we describe a method for constructing knots that can be unknotted via a single twist, and whose $d$-fold covers admit $SU(2)$-representations.

\begin{definition}\label{defn: braid adapted to a SICUP matrix}
 A tangle $\beta$ is said to be {\bf adapted} to a $d\times d$ SICUP matrix $A=(a_{ij})_{d\times d}$ if the following hold.
\begin{enumerate}
	\item The closure $K = \widehat{\beta}$ is the unknot.
	\item The closure $L = \widehat{\beta^d}$ is a $d$-component link with linking matrix $A$.
	\item If $L_1$ is the first component of $L$, then either
            \begin{itemize}
            \item $\nu^\sharp(L_i)\neq 0$ and $a_{ii}\leq \nu^\sharp(L_i)$, or
            \item $\nu^\sharp(L_i)=0$, $L_i$ is $V$-shaped, and $a_{ii}\leq \nu^\sharp(L_i)+1$, or
            \item $\nu^\sharp(L_i)=0$, $L_i$ is $W$-shaped and $a_{ii}\leq \nu^\sharp(L_i)-1$.
            \end{itemize}
\end{enumerate}
\end{definition}

%
%

\begin{remark}\label{neg-to-pos} In \Cref{defn: sicup matrices}, the matrix $A$ is stipulated to be positive definite, whereas in \Cref{thm:nu} the linking matrix must be negative definite. This mismatch is resolved by taking the mirror of the closure of a tangle adapted to an SICUP matrix: mirroring reverses every surgery slope and every linking number, replacing $A$ with $-A$ and thus transforming a positive-definite matrix into a negative-definite one. Furthermore, as shown in \cite{baldwin2020concordance}, $\nu^{\sharp}(\overline{K}) = -\nu^{\sharp}(K)$ while the knot’s shape is preserved. Hence, after taking the mirror image, all hypotheses of \Cref{thm:nu} will be satisfied.
\end{remark}

\Crefparts{mainthm}{nu} is precisely the following statement:
\begin{proposition}\label{existence_sicup}
Suppose $A$ is a $d\times d$ SICUP matrix and $\gamma$ is a tangle adapted to $A$. Let $\alpha$ be the braid axis for $\gamma$ (as in the left of \Cref{fig: exmp of 10 braid}) 
and let $K$ be the image of $\alpha$ in $S^3_{+1}(\widehat{\gamma})\cong S^3$. Then the $d$-fold branched cover $\Sigma=\Sigma_d(K)$ is an integer homology sphere, and its fundamental group $\pi_1(\Sigma)$ admits an irreducible $SU(2)$-representation.
\end{proposition}

\begin{proof} The hypotheses imply that $K$ can be unknotted by a single twist along $\gamma$ (where $\gamma$ here denotes both the braid and its closure by a slight abuse of notation). Then, $\gamma$ being adapted to $A$ means that $A$ is the linking matrix of the equivariant surgery description of $\Sigma$ obtained from lifting $\gamma$ as in \Cref{surgery_cover}. Since $A$ is unimodular and $|H_1(\Sigma;\Z)|=|\det(A)|$, $\Sigma$ is an integer homology sphere. The existence of irreducible representations follows from \Cref{thm:nu} (modulo the changes explained in \Cref{neg-to-pos}).
\end{proof}

\begin{exmp}
A solution to the Pell's equation $a^2-5b^2 = 4$ is $a=7$, $b=-3$. As a result, by Proposition \ref{prop: SICUP and Pell's eq}, we can construct a $5\times 5$ SICUP matrix with $x=3$, $l=-2$, and $m=1$. The braid $\gamma=[-2,-2,-2,-2,-1,2,2,2,2,2,-3,-4,5,6,-7,-8,9]$ shown on the left of \Cref{fig: exmp of 10 braid} is adapted to the SICUP matrix $A$ with first row $(3, -2, 1, 1, -2)$. Indeed, one can check by hand that the braid closure is unknotted and that the linking matrix for the closure of $\gamma^5$ is $A$. After taking a $5$-fold covering, the knot $L_1$ as in \Cref{defn: braid adapted to a SICUP matrix} is the torus knot $T(2,5)$, which is well-known to have lens-space surgeries, and hence is an L-space knot. By \cite[Theorem 1.18]{baldwin2020concordance}, we have $\nu^{\sharp}(L_1) = 3\geq x$. 

\begin{figure}[h]
\centering
\def\svgwidth{0.2875\textwidth}
\begingroup%
  \makeatletter%
  \providecommand\color[2][]{%
    \errmessage{(Inkscape) Color is used for the text in Inkscape, but the package 'color.sty' is not loaded}%
    \renewcommand\color[2][]{}%
  }%
  \providecommand\transparent[1]{%
    \errmessage{(Inkscape) Transparency is used (non-zero) for the text in Inkscape, but the package 'transparent.sty' is not loaded}%
    \renewcommand\transparent[1]{}%
  }%
  \providecommand\rotatebox[2]{#2}%
  \newcommand*\fsize{\dimexpr\f@size pt\relax}%
  \newcommand*\lineheight[1]{\fontsize{\fsize}{#1\fsize}\selectfont}%
  \ifx\svgwidth\undefined%
    \setlength{\unitlength}{205.64440918bp}%
    \ifx\svgscale\undefined%
      \relax%
    \else%
      \setlength{\unitlength}{\unitlength * \real{\svgscale}}%
    \fi%
  \else%
    \setlength{\unitlength}{\svgwidth}%
  \fi%
  \global\let\svgwidth\undefined%
  \global\let\svgscale\undefined%
  \makeatother%
  \begin{picture}(1,1.51890346)%
    \lineheight{1}%
    \setlength\tabcolsep{0pt}%
    \put(0,0){\includegraphics[width=\unitlength,page=1]{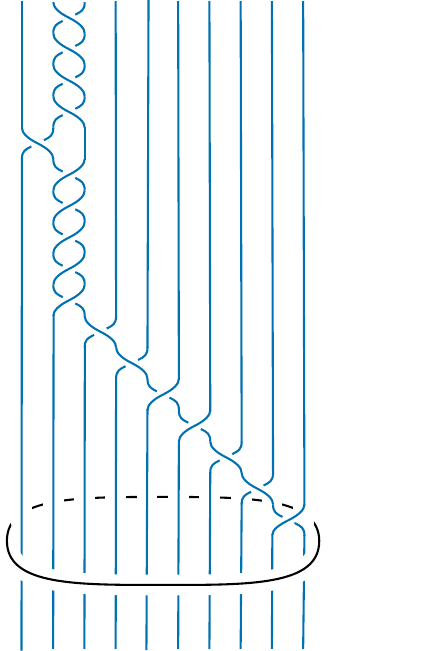}}%
    \put(0.75661492,0.47114716){\color[rgb]{0,0.44705882,0.69803922}\makebox(0,0)[lt]{\lineheight{1.25}\smash{\begin{tabular}[t]{l}$\gamma$\end{tabular}}}}%
    \put(0.73714161,0.91481795){\color[rgb]{0,0.44705882,0.69803922}\makebox(0,0)[lt]{\lineheight{1.25}\smash{\begin{tabular}[t]{l}$+1$\end{tabular}}}}%
    \put(0.75661492,0.26202231){\color[rgb]{0,0,0}\makebox(0,0)[lt]{\lineheight{1.25}\smash{\begin{tabular}[t]{l}$K$\end{tabular}}}}%
  \end{picture}%
\endgroup%
\quad
\def\svgwidth{0.2875\textwidth}
\begingroup%
  \makeatletter%
  \providecommand\color[2][]{%
    \errmessage{(Inkscape) Color is used for the text in Inkscape, but the package 'color.sty' is not loaded}%
    \renewcommand\color[2][]{}%
  }%
  \providecommand\transparent[1]{%
    \errmessage{(Inkscape) Transparency is used (non-zero) for the text in Inkscape, but the package 'transparent.sty' is not loaded}%
    \renewcommand\transparent[1]{}%
  }%
  \providecommand\rotatebox[2]{#2}%
  \newcommand*\fsize{\dimexpr\f@size pt\relax}%
  \newcommand*\lineheight[1]{\fontsize{\fsize}{#1\fsize}\selectfont}%
  \ifx\svgwidth\undefined%
    \setlength{\unitlength}{210bp}%
    \ifx\svgscale\undefined%
      \relax%
    \else%
      \setlength{\unitlength}{\unitlength * \real{\svgscale}}%
    \fi%
  \else%
    \setlength{\unitlength}{\svgwidth}%
  \fi%
  \global\let\svgwidth\undefined%
  \global\let\svgscale\undefined%
  \makeatother%
  \begin{picture}(1,1.5)%
    \lineheight{1}%
    \setlength\tabcolsep{0pt}%
    \put(0,0){\includegraphics[width=\unitlength,page=1]{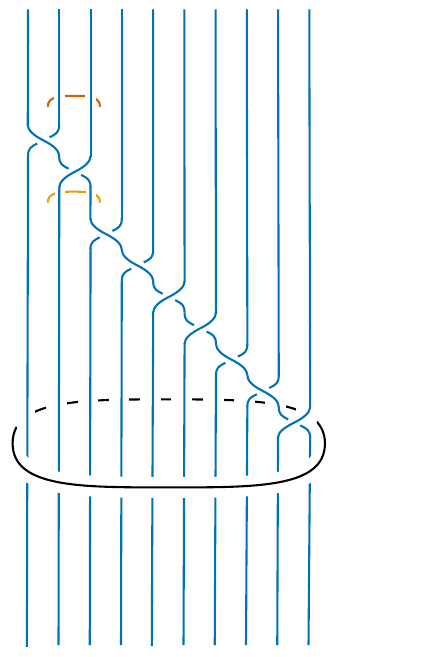}}%
    \put(0.74862249,0.63710574){\color[rgb]{0,0.44705882,0.69803922}\makebox(0,0)[lt]{\lineheight{1.25}\smash{\begin{tabular}[t]{l}$\gamma$\end{tabular}}}}%
    \put(0.72955308,0.92904085){\color[rgb]{0,0.44705882,0.69803922}\makebox(0,0)[lt]{\lineheight{1.25}\smash{\begin{tabular}[t]{l}$+1$\end{tabular}}}}%
    \put(0,0){\includegraphics[width=\unitlength,page=2]{SICUP-annulus.pdf}}%
    \put(0.25389764,1.24586985){\color[rgb]{0.83529412,0.36862745,0}\makebox(0,0)[lt]{\lineheight{1.25}\smash{\begin{tabular}[t]{l}$+1/2$\end{tabular}}}}%
    \put(0,0){\includegraphics[width=\unitlength,page=3]{SICUP-annulus.pdf}}%
    \put(0.26311306,1.015779){\color[rgb]{0.90196078,0.62352941,0}\makebox(0,0)[lt]{\lineheight{1.25}\smash{\begin{tabular}[t]{l}$-1/2$\end{tabular}}}}%
    \put(0.74862249,0.50407046){\color[rgb]{0,0,0}\makebox(0,0)[lt]{\lineheight{1.25}\smash{\begin{tabular}[t]{l}$K$\end{tabular}}}}%
    \put(0,0){\includegraphics[width=\unitlength,page=4]{SICUP-annulus.pdf}}%
  \end{picture}%
\endgroup%
\quad
\def\svgwidth{0.35\textwidth}
\begingroup%
  \makeatletter%
  \providecommand\color[2][]{%
    \errmessage{(Inkscape) Color is used for the text in Inkscape, but the package 'color.sty' is not loaded}%
    \renewcommand\color[2][]{}%
  }%
  \providecommand\transparent[1]{%
    \errmessage{(Inkscape) Transparency is used (non-zero) for the text in Inkscape, but the package 'transparent.sty' is not loaded}%
    \renewcommand\transparent[1]{}%
  }%
  \providecommand\rotatebox[2]{#2}%
  \newcommand*\fsize{\dimexpr\f@size pt\relax}%
  \newcommand*\lineheight[1]{\fontsize{\fsize}{#1\fsize}\selectfont}%
  \ifx\svgwidth\undefined%
    \setlength{\unitlength}{260bp}%
    \ifx\svgscale\undefined%
      \relax%
    \else%
      \setlength{\unitlength}{\unitlength * \real{\svgscale}}%
    \fi%
  \else%
    \setlength{\unitlength}{\svgwidth}%
  \fi%
  \global\let\svgwidth\undefined%
  \global\let\svgscale\undefined%
  \makeatother%
  \begin{picture}(1,1.21153846)%
    \lineheight{1}%
    \setlength\tabcolsep{0pt}%
    \put(0,0){\includegraphics[width=\unitlength,page=1]{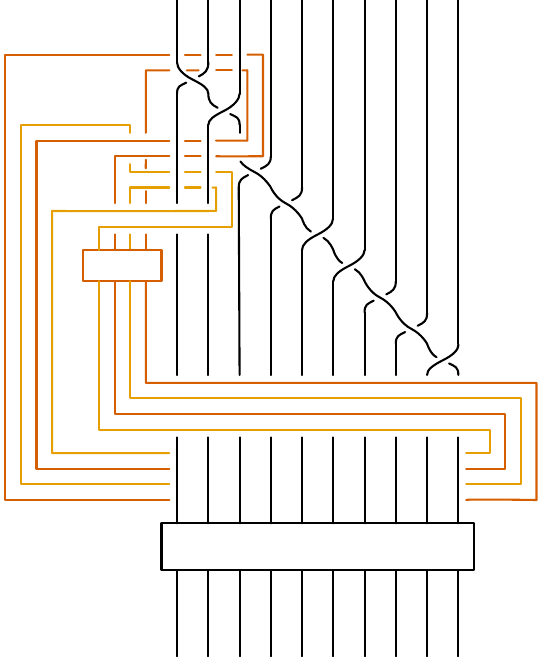}}%
    \put(0.05256808,1.13124586){\color[rgb]{0.83529412,0.36862745,0}\makebox(0,0)[lt]{\lineheight{1.25}\smash{\begin{tabular}[t]{l}$+1/2$\end{tabular}}}}%
    \put(0.05094876,0.99422662){\color[rgb]{0.90196078,0.62352941,0}\makebox(0,0)[lt]{\lineheight{1.25}\smash{\begin{tabular}[t]{l}$-1/2$\end{tabular}}}}%
    \put(0.54756611,0.17841046){\color[rgb]{0,0,0}\makebox(0,0)[lt]{\lineheight{1.25}\smash{\begin{tabular}[t]{l}$-1$\end{tabular}}}}%
    \put(0.17591271,0.70442458){\color[rgb]{0.83529412,0.36862745,0}\makebox(0,0)[lt]{\lineheight{1.25}\smash{\begin{tabular}[t]{l}$-1$\end{tabular}}}}%
    \put(0.85429687,0.86975661){\color[rgb]{0,0,0}\makebox(0,0)[lt]{\lineheight{1.25}\smash{\begin{tabular}[t]{l}$K$\end{tabular}}}}%
  \end{picture}%
\endgroup%

\caption{Left: A 10-braid $\gamma$ adapted to the SICUP with first row $(3,-2,1,1,-2)$, and a knot $K$ with untwisting sequence given by the $+1$-framed $\widehat{\gamma}$. Center: The knot $K$ presented in terms of a $2$-annulus twist and an unknotting sequence. Right: The knot $K$ whose $5$-fold branched cover admits an irreducible $SU(2)$ by \Cref{existence_sicup}.}\label{fig: exmp of 10 braid}
\end{figure}
\end{exmp}

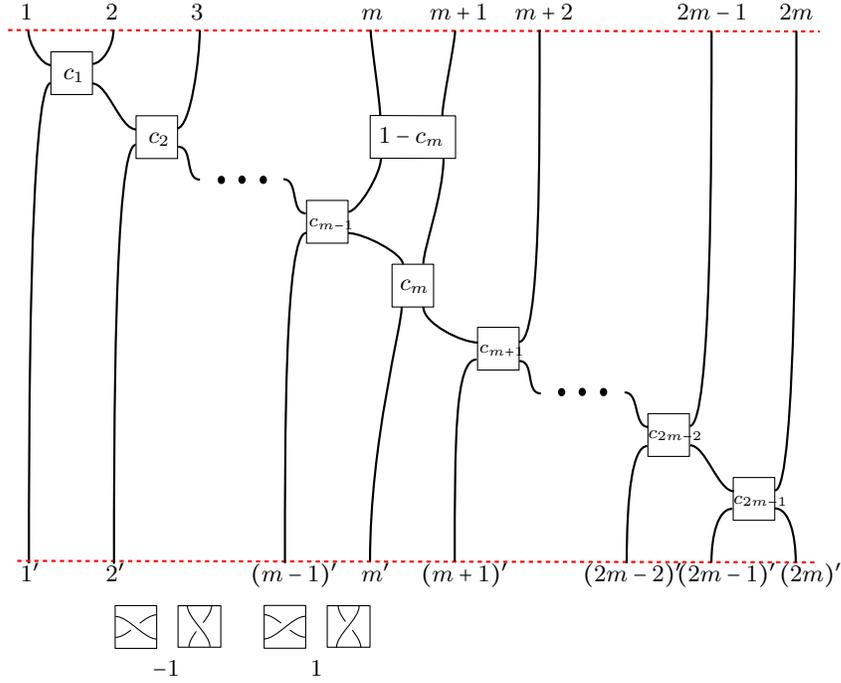
\begin{figure}[h!]
\fontsize{9}{10}\selectfont
\centering
	\begin{overpic}[width = 0.75\textwidth]{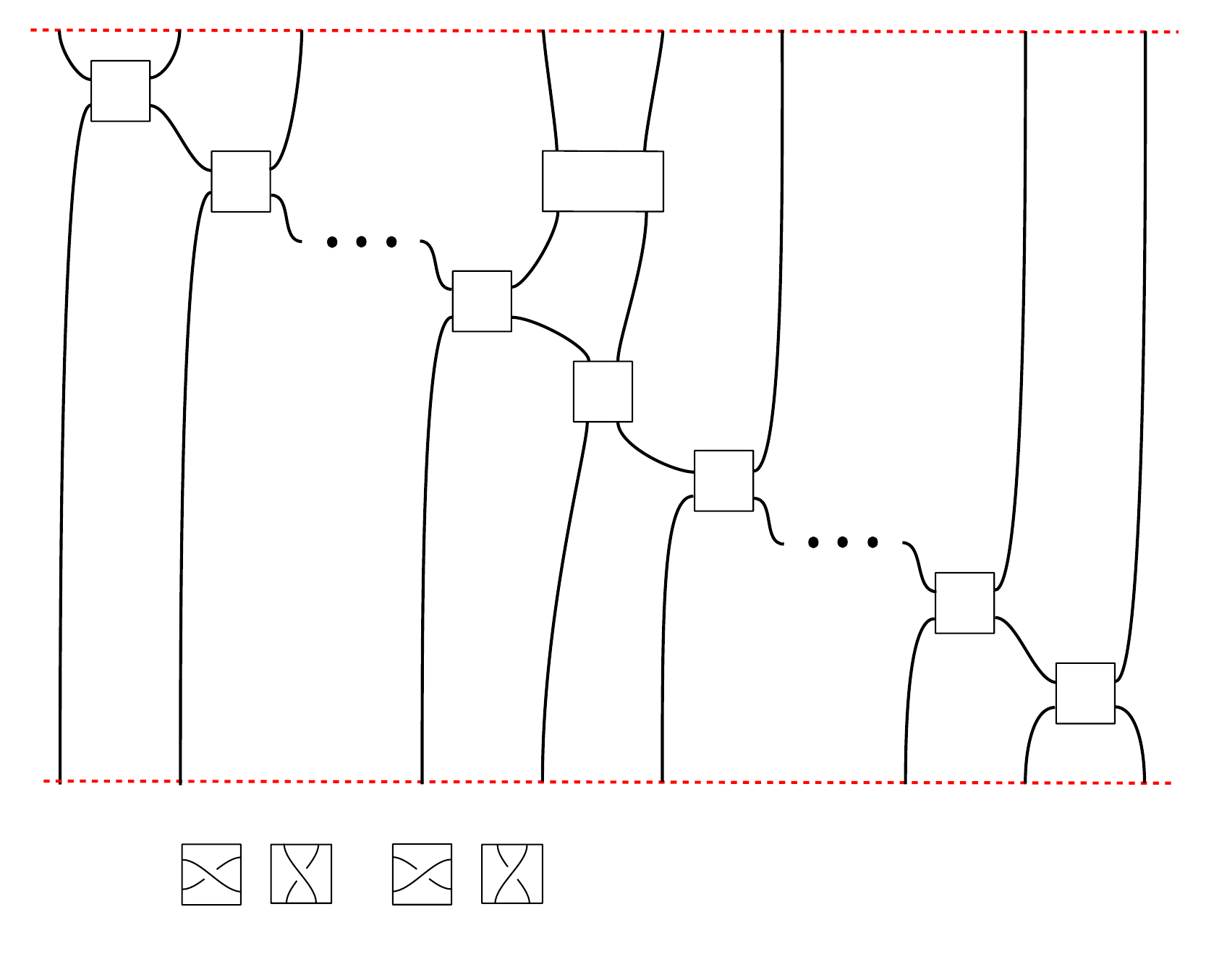}
		\put(19.5,3){$-1$}
		\put(38,3){$1$}
		
		\put(9,73){$c_1$}
		\put(19,65.5){$c_2$}
		\put(37.8,55.7){\tiny $c_{m-1}$}
		\put(48.5,48.2){$c_m$}
		\put(46,65.5){$1-c_m$}
		\put(57.8,40.9){\tiny $c_{m+1}$}
		\put(77.6,30.9){\tiny $c_{2m-2}$}
		\put(87.6,23.2){\tiny $c_{2m-1}$}
		
		\put(4,80){$1$}
		\put(14,80){$2$}
		\put(24,80){$3$}
		\put(44,80){$m$}
		\put(52,80){$m+1$}
		\put(62,80){$m+2$}
		\put(81,80){$2m-1$}
		\put(93,80){$2m$}
		
		\put(4,14){$1'$}
		\put(14,14){$2'$}
		\put(31,14){$(m-1)'$}
		\put(44,14){$m'$}
		\put(51,14){$(m+1)'$}
		\put(70,14){$(2m-2)'$}
		\put(81,14){$(2m-1)'$}
		\put(93,14){$(2m)'$}
	\end{overpic}
	\caption{The tangle $\sigma = \sigma(c_1,\cdots,c_{2m-1})$. The $c_i$'s represent odd numbers. Notice that the boxes with labels $c_m$ and $1-c_m$ represent `vertical' twists, whereas the rest represent `horizontal' ones.}\label{fig: the tangle}
\end{figure}

\begin{proposition}
Suppose $m$ is a positive integer and $\sigma$ is the $2m$-strand tangle
\[\sigma = \sigma(c_1,\cdots,c_{2m-1})\]
as in Figure \ref{fig: the tangle}. Then $\sigma$ is adapted to an $m\times m$ SICUP matrix $A=(a_{i,j})_{m\times m}$ if the following statements hold:
\begin{itemize}
	\item For $j\in\{1,\cdots,2m-1\}$, $c_j$ is an odd integer.
	\item $a_{11}\leq c_m -2$.
	\item $a_{12}=
	\frac{1}{2}\left(c_{1}+c_{m+1}+c_{m-1}+c_{2m-1}-c_m+1\right)$
	\item $a_{1j}=\frac{1}{2}\left(c_{j-1}+c_{m+j-1}+c_{m-j+1}+c_{2m-j+1}\right)$ for $j\in\{3,\cdots,\lceil\frac{m+1}{2}\rceil\}$
\end{itemize}
\end{proposition}

\begin{figure}[h]
\fontsize{9}{10}\selectfont
\centering
\def\svgwidth{0.325\textwidth}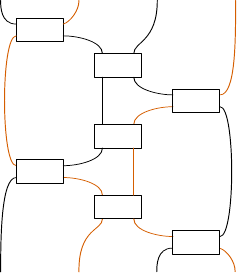\qquad\qquad\qquad
\def\svgwidth{0.325\textwidth}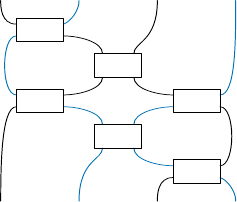
\caption{Left: The tangle representation for $L_1\cup L_2$. The black strands represent $L_1$ and the red strands represent $L_2$. Right: The tangle representation for $L_1\cup L_j$. The black strands represent $L_1$ and the blue strands represent $L_j$..}\label{fig:L1_L2}
\end{figure}

\begin{proof}
Note that since all $c_i$ are odd, in Figure \ref{fig: the tangle}, the end $1$ at top is connected to $(2m)'$ at bottom, and for all $j\in\{2,\cdots,2m\}$ the end $j$ at top is connected to the end $(j-1)'$ at bottom.

We first show that the closure $\hat{\sigma}$ is the unknot. When closing it up, the ends $1$ and $1'$ are connected. And hence $|c_1|$ Reidemeister I moves will resolve all crossings in the block $c_1$. Then $|c_2|$ Reidemeister I moves will resolve all crossings in the block $c_2$. At the same time, the ends $2m$ and $(2m)'$ are also connected, and hence $|c_{2m-1}|$ Reidemeister I moves will resolve all crossings in the block $c_{2m-1}$. Finally, we obtained a two-strand braid on $m$ and $(m+1)$, and the twist between them is $1-c_m + c_m = 1$. Hence $\hat{\sigma}=T(2,1)$, which is the unknot.

Second, let
\[
L = \widehat{\sigma^m} = L_1\cup\cdots\cup L_m.
\]
Tangle representations for $L_1\cup L_2$ and $L_1\cup L_j$ $(j>2)$ are included in \Cref{fig: L1_L2}. 
Ignoring the non-black strands in those gives the braid representation $(2,c_m)$ for $L_1$, thus showing that $L_1=T(2,c_m)$.  This is an instanton $L$-space knot and according to the definition of $\nu^{\sharp}$ in \cite{baldwin2020concordance} and \cite[Theorem 1.1]{lidman2022framed}, we conclude
\[
\nu^{\sharp}(L_1) = 2g(L_1) - 1 = c_m - 2\geq a_{11}.
\]
%
%
%

The conditions on the linking numbers can be verified using \Cref{fig:L1_L2}
Using those, it is straightforward to compute that
\begin{align*}
{\rm lk}(L_1,L_2) &= \frac{c_{1}+c_{m+1}+c_{m-1}+c_{2m-1} + 1-c_m}{2}, \text{ and }\\
{\rm lk}(L_1,L_j) &= \frac{c_{j-1}+c_{m+j-1}+c_{m-j+1}+c_{2m-j+1}}{2} \text{ for } j\neq 2,m.
\end{align*}

\end{proof}
%
%
%
%
%
%

\bibliographystyle{alpha}
\bibliography{covers-reps-references.bib}

\end{document}